\documentclass[oneside,10pt]{article}          % please do not change
\usepackage[b5paper]{geometry}      % your paper can be easily printed on a4 or letter paper with enlargenment
\usepackage{amsfonts,amsmath,latexsym,amssymb} % these packages are required
\usepackage{amsthm}                % please use one of this two options for theorems
\usepackage{mathrsfs,upref}         % not so essential, but part of journal style
\usepackage{mathptmx}           % Journal is printed with poscript fonts:
\usepackage{indentfirst}
\newtheorem{lemma}{Lemma}[section]

\newtheorem{theorem}{Theorem}[section]         % alternative set of theorem numbering

\theoremstyle{definition}

\newtheorem{remark}{Remark}

\numberwithin{equation}{section}
\begin{document}

\title{
Schauder estimates for solutions of sub-Laplace equations with Dini terms
\thanks{This work was supported by the National
Natural Science Foundation of China (Grant No. 11271299), Natural
Science Foundation Research Project of Shaanxi Province (Grant No.
2012JM1014).} }

% Short title is optional, it will appear in running heads.
% It is necessary only if the title is to long to be used in running heads

\author{Tingxi Hu, Pengcheng Niu\thanks{ Corresponding author. pengchengniu@nwpu.edu.cn(P. Niu)}\\Department of Applied Mathematics,\\Northwestern Polytechnical University, Xi'an 710129, China}

%\address{Third Author, Full postal Address of the Third Author\\
%\email{email of the Third Author}}

%\dedicated{Dedicated to...}                    % Optional

%\date{07.03.2013}                               % Please, write the date of submission

        %% AMS-2010 subj class. The list can be found on http://www.ams.org/mathscinet/msc/msc2010.html

%%\thanks{This research is supported by ...\par This paper is a lecture that was given at...}
        %% Optional. Only one command thanks is allowed, use \par inside text if you need multiple thanks.
\maketitle
\begin{abstract}
In this paper we establish Schauder estimates for the sublalpace equation%/
\[\Sigma _{j = 1}^mX_j^2u = f,\]
where ${X_1},{X_2}, \ldots ,{X_m}$ is a system of smooth vector field which generates the first
layer in the Lie algebra of a Carnot group. We drive the estimate for the second
order derivatives of the solution to the equation with Dini continue inhomogeneous term $f$ by the perturbation argument.

\textbf{Keywords:}\;Carnot group; sub-Laplace; Schauder estimate; Dini continue; perturbation argument.

\textbf{MSC2010:}\;35B65; 35R03.
\end{abstract}

%%%%% END OF TITLE PAGE %%%%%%%%%%%%%%%%%%%%%%%%%%%%%%%%%%%%%%%%%%%%%%

%%%%% BODY OF THE PAPER %%%%%%%%%%%%%%%%%%%%%%%%%%%%%%%%%%%%%%%%%%%%%%
% You should eventually delete (after reading) the rest of the text below %%

\section{Introduction}
Schauder estimates play an important role in the theory of elliptic equations, see [6, 11]. For the second
order uniformly elliptic equation in any bounded domain $\Omega  \subset {\mathbb{R}^n}$%/
\[\Sigma _{i,j = 1}^n{a_{ij}}(x)\partial _{ij}^2u = f,\]
such estimates provide a bound of the H\"{o}lder norm in $\Omega$ of the second derivatives of the solution
$u$ in terms of the H\"{o}lder norms in $\Omega$ of the coefficients $a_{ij}$ and $f$.

A sharper form of these estimates was introduced by Caffarelli [3]
in the study of fully non-linear elliptic equations. He derived
Schauder estimates for viscosity solutions by comparing the
solutions with osculating quadratic polynomials in a neighborhood of
a fixed point, and this method is called \lq\lq perturbation argument". In
Caffarelli's approach, the H\"{o}lder regularity of $u$ at a point is
basically determined by the H\"{o}lder regularity of $a_{ij}$ and
$f$ at the same point, hence such estimates are said pointwise
Schauder estimates. Wang in [23] compared the quadratic part of the
solutions to the Laplace equation with solutions of approximate
equations and proved the H\"{o}lder norm of $D^2u$ in terms of
Dini continuous inhomogeneous term. Afterwords, the method was
used to investigate the fully non-linear elliptic and parabolic
equations, see Liu-Trudinger-Wang[16], Tian-Wang [22].

For degenerate elliptic equations constructed by left translation invariant vector fields, several authors derived
Schauder estimates, see Lunardi in [19], Capogna-Han [5], Polidoro
-Di Francesco [21] and Guti\'{e}rrez-Lanconelli [12].
Schauder estimates for heat type euqations induced by smooth vector fields satisfying H\"{o}rmander's finite rank condition were
showed by Bramanti-Brandolini [2].

Recently, Jiang-Tian [15] showed Schauder estimates for the Kohn-Laplace equation with Dini continuous inhomogeneous
term in the Heisenberg group in the spirit of [23]. We generalize the result in [15] to the sub-Laplace equations in
Carnot groups.

In the present paper we consider the equation
\begin{equation}
Lu \equiv \Sigma _{j = 1}^mX_j^2u = f  \;\;\;in \;{B_1}(0),
\end{equation}
in which $L$ is the sub-laplacian on a Carnot group $G$ and the right hand term $f$ is Dini continuous, i.e., $f$ satisfies%/
\[\int_0^1 {\frac{{{\omega _f}(r)}}{r}dr}  < \infty ,\]
where ${\omega _f}(r) = \mathop {\sup }\limits_{d(\xi ,\eta ) < r} |f(\xi ) - f(\eta )|$ and $d(\xi ,\eta )$ is the
pseudo-distance (see next section) between $\xi$ and $\eta$, ${B_1}(0)$ denotes the unit gauge ball centered at origin.

Our main result is the following
\begin{theorem}
Let $u \in {C^2}({B_1}(0))$ be a solution of (1.1), then for any $\xi ,\eta  \in {B_{1/2}}(0)$, $d = d(\xi ,\eta )$,
there exists a positive constant $C$ such that%/
\begin{eqnarray}
&&|{X_i}{X_j}u(\xi ) - {X_i}{X_j}u(\eta )| \nonumber\\
&&\leqslant C\left(d(\mathop {\sup }\limits_{{B_1}(0)} |u| +
||f|{|_{{L^\infty }}} + \int_{\sqrt d }^1 {\frac{{{\omega
_f}(r)}}{{{r^2}}}dr} ) + \int_0^{\sqrt d } {\frac{{{\omega
_f}(r)}}{r}dr}\right).
\end{eqnarray}

In particular, if $f \in {C^{0,\alpha }}({B_1}(0))\;(0 < \alpha \leqslant 1)$, then
\begin{equation}
|{X_i}{X_j}u(\xi ) - {X_i}{X_j}u(\eta )| \leqslant C{d^{\alpha
/2}}\left(\mathop {\sup }\limits_{{B_1}(0)} |u| +
||f|{|_{{C^{0,\alpha }}}}\right), \;\; \alpha  \in (0,1),
\end{equation}
\begin{eqnarray}
&&|{X_i}{X_j}u(\xi ) - {X_i}{X_j}u(\eta )| \nonumber\\
&&\leqslant C{d^{1/2}}\left(\mathop {\sup }\limits_{{B_1}(0)} |u| +
||f|{|_{{C^{0,1}}}}\left(1 + |\sqrt d \log \sqrt d |\right)\right),
\;\;\alpha=1.
\end{eqnarray}
\end{theorem}

The plan of the paper is as follows: in Section 2 we introduce knowledge related to Carnot groups and some
preliminary lemmas. Also a maximum principle for (1.1) with Dirichlet boundary value problem is proved.
Section 3 is devoted to the proof of Theorem 1.1. We mention that the treatment for the Taylor polynomials
in the Carnot group is more complicated than in the Heisenberg group. Also necessary techniques to use the
perturbation argument are given.

\section{Preliminary results}
We begin by describing several known facts on Carnot groups and refer to [1, 10] for more information. Especially,
we provide a maximum principle (Lemma 2.5) for solutions of a boundary value problem to the sub-Laplace equation.

A Carnot group $G$ of step $s$ is a simply connected nilpotent Lie group such that its Lie algebra $\mathfrak{g}$
admits a stratification $\mathfrak{g}= \oplus _{l = 1}^s{V_l}$, with $[{V_1},{V_l}] = {V_{l + 1}}$ $(l = 1,2, \ldots ,s - 1)$
and $[{V_1},{V_s}] = \left\{ 0 \right\}$. Denoting ${m_l} = \dim {V_l}$, we fix on $G$ a system of coordinates
$\xi  = \left( {{z_1},{z_2}, \ldots ,{z_s}} \right)$, in which ${z_l} = ({x_{l,1}},{x_{l,2}}, \ldots ,{x_{l,{m_l}}}) \in {\mathbb{R}^{{m_l}}}$.

Every Carnot group $G$ is naturally equipped with a family of non-isotropic dilations defined by $\delta_r$:%/
$${\delta _r}(\xi ) = (r{z_1},{r^2}{z_2}, \ldots ,{r^s}{z_s})\;,\xi  \in G,r > 0,$$
and the homogeneous dimension of $G$ is given by $Q = \sum\limits_{l
= 1}^s {l{m_l}}$. We express by $dH(\xi)$ a fixed bi-invariant Haar
measure on $G$. One easily sees $dH({\delta _r}(\xi )) = {r^Q}dH(\xi
)$. By the Baker-Campbell-Hausdorff formula, the group law on $G$  is
$$\xi \eta  = \xi  + \eta  + \sum\limits_{1 \leqslant l,k \leqslant s} {{Z_{l,k}}(\xi ,\eta )} \;,\;\;\xi ,\eta  \in
G,
$$
where ${Z_{l,k}}(\xi ,\eta )$  is a fixed linear combination of
iterated commutators containing $l$ times $\xi$  and $k$  times
$\eta$.

The homogenous norm of $\xi$  on $G$  is defined by  $|\xi | =
{(\Sigma _{j = 1}^s|{z_j}{|^{2s!/j}})^{1/2s!}}$, where $|{z_j}|$
denotes the Euclidean norm of  ${z_j} \in {\mathbb{R}^{{m_j}}}$. Such
homogenous norm on  $G$ can be used to define a pseudo-distance on
$G$ which is  $d(\xi ,\eta ) = |{\xi ^{ - 1}}\eta |$. Denote the
gauge ball of radial  $r$ centered at  $\xi$ by ${B_r}(\xi ) = \{
\eta  \in G|d(\xi ,\eta ) < r\}$.

Let $X = \{ {X_1},{X_2}, \ldots ,{X_m}\}$  be a basis of  $V_{1}$,
then we can write $X_{i}$  as
$${X_i} = {\partial _{1,i}} + \sum\limits_{j= i+1}^m{a_{ij}}(\xi ){\partial _{1,j}}
+\sum\limits_{l =2}^s\sum\limits_{k= 1}^{m_l}{b_{ilk}}(\xi ){\partial _{l,k}},\;{X_i}(0) = {\partial
_{1,i}},
$$
where ${a_{ij}}(\xi )$ and ${b_{ilk}}(\xi )$ are polynomials. The sub-Laplacian on $G$
associated with  $X$ is of the form $L = \Sigma _{i = 1}^mX_i^2$,
which is a hypoelliptic second order partial differential operator.
Obviously, it holds  $L(\frac{{x_{1,1}^2}} {2}) = 1$.

Suppose that $\{ {X_{l,1}},{X_{l,2}}, \ldots ,{X_{l,{m_l}}}\}$ is a
basis of  $V_{l}$ and consider a multi-index $I = \{ ({i_k},{j_k})\}
_{k = 1}^s$ with  ${i_k} \in \{ 1,2, \ldots ,l\} ,{j_k} \in \{ 1,2,
\ldots ,{m_{{i_k}}}\}$. For a smooth function  $f$ on  $G$, we
denote a derivative of  $f$ with order  $|I| = \Sigma _{k = 1}^s{i_k}$ by
$${X^I}f = {X_{{i_1},{j_1}}}{X_{{i_2},{j_2}}} \ldots {X_{{i_s},{j_s}}}f.
$$
We will use the notation $Xf = ({X_1}f,{X_2}f, \ldots ,{X_m})$ for
convenience too.

A polynomial on  $G$ is a function which can be expressed in the
exponential coordinates by
$$P(\xi ) = \sum\limits_J {{a_J}{x^J}},$$
where  $J = \{ {j_{i,k}}\} _{i = 1, \ldots ,{m_k}}^{k = 1, \ldots
,s}$,  ${a_J}$ are the real numbers and  ${x^J} = \Pi _{i = 1,
\ldots ,{m_k}}^{k = 1, \ldots ,s}x_{i,k}^{{j_{i,k}}}$. The
homogeneous degree of the monomial  $x^{J}$ is given by the sum $|J|
= \Sigma _{k = 1}^s\Sigma _{i = 1}^{{m_k}}k{j_{i,k}}$.

Let $\Omega  \subset G$  be an open set. If   $k \in N$ and  $1
\leqslant p < \infty$, we define the horizontal Sobolev space by
$$H{W^{k,p}}(\Omega ) = \{ f:\;|{X^I}f| \in {L^p}(\Omega ),0 \leqslant |I| \leqslant k\}
.$$ Then we illustrate the H\"{o}lder space and Lipschitz space with
respect to the pseudo-distance. If $0 < \alpha  \leqslant 1$ and $f$
is a function defined in an open set  $\Omega$, let
$${[f]_{{C^{0,\alpha }}}} = \sup \left\{ {\frac{{|f(\xi ) - f(\eta )|}}
{{d{{(\xi ,\eta )}^\alpha }}}:\xi ,\eta  \in \Omega ,\xi  \ne \eta }
\right\}.
$$
The H\"{o}lder space is defined by
$${C^{0,\alpha }}(\Omega ) = \{ f|\;[f]_{{{C^{0,\alpha }}}} < \infty
\}, \;0<\alpha<1 $$
and Lipschitz space by ${C^{0,1}}(\Omega ) = \{
f|\;[f]_{{{C^{0,1}}}} < \infty \}$. In addition we denote that
$||f||_{{{C^{0,\alpha }}}}:=[f]_{{{C^{0,\alpha }}}}+ ||f||_{L^\infty}$, for $0< \alpha \leqslant 1$.

We introduce some known results that will be used in this paper.
\begin{lemma}([1, pp.390-391])
The gauge balls  ${B_r}(\xi )(\xi  \in G,r > 0)$
are  $L$-regular open sets, i.e., for any  $f \in
{C^\infty }({B_r}(\xi ))$, there exists a Perron-Wiener-Brelot
generalized solution  $u \in {C^\infty }({B_r}(\xi )) \cap
C(\overline {{B_r}(\xi )} )$ to the boundary value problem
\begin{equation}\left\{ {\begin{array}{*{20}{c}}
   {Lu = f}  \\
   {u{|_{\partial {B_r}(\xi )}} = g}  \\

 \end{array} } \right.\;\;\begin{array}{*{20}{c}}
   {in\;{B_r}(\xi ),}  \\
   {g \in C(\partial {B_r}(\xi )).}  \\

 \end{array}
\end{equation}
\end{lemma}

\begin{lemma}(a priori estimates, [4])
Let $\Omega\subset G$  be an open set and $u$  be $L$  harmonic,
i.e.,  $u$ satisfies  $Lu=0$, then for a given integer  $k$ and any
multiple-index $I$, $|I| \leqslant k$, there exists a constant $C$
depending on  $G$ and $k$ such that if  $\overline {{B_r}(\xi )}
\subset \Omega $, then
\begin{equation}
|{X^I}u|(\eta ) \leqslant C{r^{ - k}}\mathop {\sup
}\limits_{\overline {{B_r}(\xi )} } |u|,\;\;\eta \in {B_r}(\xi ).
\end{equation}
\end{lemma}

\begin{lemma} (Folland-Stein [7])
Let $\Omega\subset G$ be an open set, then for any  $1<p<Q$, there
exist a positive constant $S_{p}$ depending on  $G$, such that for
$f \in C_0^\infty (\Omega )$,
\begin{equation}
{\left(\int_\Omega  {|f{|^{p*}}dH} \right)^{1/p*}} \leqslant
{S_p}{\left(\int_\Omega {|Xf{|^p}dH} \right)^{1/p}},
\end{equation}where $p* = \frac{{pQ}} {{Q - p}}, |Xf| = {(\Sigma _{j =
1}^m|{X_j}f{|^2})^{1/2}}.$

\end{lemma}

The following technical lemma is adapted from Chen and Wu [6].
\begin{lemma}(De Giorgi's iteration lemma, [6])
Let $\varphi (t)$  be a nonnegative and non-increasing function on
$[{k_0}, + \infty )$ satisfying
$$\varphi (h) \leqslant \frac{C}
{{{{(h - k)}^\alpha }}}{[\varphi (k)]^\beta }, \;h > k \geqslant
{k_0},
$$
for some constant  $C > 0,\alpha  > 0,\beta  > 1$. Then we have
\begin{equation}
\varphi ({k_0} + \tilde d) = 0,
\end{equation}
in which $\tilde d = {C^{1/\alpha }}{[\varphi ({k_0})]^{(\beta  -
1)/\alpha }}{2^{\beta /(\beta  - 1)}}$.
\end{lemma}
Following the method of proving a classical maximum principle in [6,
Theorem 2.4], we can obtain the following result by combining Lemmas
2.3 and 2.4.
\begin{lemma}  (Maximum principle)
Let $\Omega\subset G$  be an open set, $f \in {L^\infty }(\Omega )$,
 $u \in {C^2}(\Omega )$ solves (2.1), then
\begin{equation}
\mathop {\sup }\limits_\Omega  |u| \leqslant \mathop {\sup
}\limits_{\partial \Omega } |g| + C||f|{|_{{L^\infty }(\Omega
)}}|\Omega {|^{2/Q}}{2^{(Q + 2)/4}}.
\end{equation}
\end{lemma}

\textbf{Proof.}  Notice that for every  $\varphi  \in C_0^2(\Omega
)$, we have
$$\int_\Omega  {\Sigma _{i = 1}^{m}{X_i}u{X_i}\varphi } dH =  - \int_\Omega  {f\varphi dH}.
$$
Set   ${k_0} = \mathop {\sup }\limits_{\partial \Omega } |g|$ and
$\varphi  = {(u - k)_ + }$ with  $k>k_{0}$, and denote  $A(k) = \{
\xi  \in \Omega |u > k\}$. It is easy to know ${X_i}\varphi  =
{X_i}u$  in  $A(k)$. Then
\begin{equation}
\int_{A(k)} {\Sigma _{i = 1}^{m}|{X_i}\varphi {|^2}dH}  =
\int_{A(k)} {\Sigma _{i = 1}^{m}{X_i}u{X_i}\varphi } dH =
\int_{A(k)} {f\varphi dH}.
\end{equation}

By Lemma 2.3, it obtains
\begin{equation}
{\left( {\int_{A(k)} {|\varphi {|^{2^*}}dH}}
\right)^{2/2^*}} \leqslant C\int_{A(k)} {\Sigma _{i =
1}^{m}|{X_i}\varphi {|^2}dH}.
\end{equation}
On the other hand,
\begin{eqnarray}
\int_{A(k)} {f\varphi dH}  &\leqslant& {\left( {\int_{A(k)}
{|\varphi {|^{2^*}}dH} } \right)^{1/2^*}}{\left( {\int_{A(k)}
{|f{|^{2Q/(2 + Q)}}dH} }
     \right)^{(2 + Q)/2Q}}
 \nonumber\\
&\leqslant& {\left(\int_{A(k)} {|\varphi {|^{2^*}}dH}
\right)^{1/2^*}}||f|{|_{{L^\infty }(\Omega )}}|A(k){|^{(2 + Q)/2Q}}.
\end{eqnarray}
Since $A(h) \subset A(k)$ and $\varphi
\geqslant h - k$ in $A(h)$ if $k<h$, it follows
\begin{equation}
{(h - k)^{2^*}}|A(h)| \leqslant \int_{A(h)} {|\varphi {|^{2^*}}dH}
\leqslant \int_{A(k)} {|\varphi {|^{2^*}}dH}.
\end{equation}
Combining (2.6)-(2.9), it yields
$$|A(h)| \leqslant \frac{{{{(C||f|{|_{{L^\infty }(\Omega )}})}^{2^*}}}}
{{{{(h - k)}^{2^*}}}}|A(k){|^{(Q + 2)/(Q - 2)}}.
$$
By Lemma 2.4 we get (2.5). $\Box$

We will need the following three Lemmas referring to [1, 8], which
are important in applying the perturbation argument.
\begin{lemma} (Taylor polynomial)
Let  $f \in {C^\infty }(G)$, then for every integer  $n$, there
exists a unique polynomial  ${P_n}(f,0)$ homogenous of degree at
most  $n$, such that
\begin{equation}
{X^I}{P_n}(f,0)(0) = {X^I}f(0),
\end{equation}
for all multiple-index  $I$ satisfying  $|I| \leqslant n$.
\end{lemma}
\begin{lemma} (Remainder in Taylor formula)
Let  $f \in {C^{n + 1}}(G), \,\xi \in G$,  then
\begin{equation}
f(\eta ) - {P_n}(f,\xi )(\eta ) = {O_{\eta  \to \xi }}({d^{n +
1}}({\xi ^{ - 1}}\eta )).
\end{equation}
\end{lemma}

\begin{lemma} (Mean value theorem)
There exist absolute constants  $b,C>0$, depending only on  $G$ and the
homogenous norm  $| \cdot |$, such that
\begin{equation}
|f(\xi \eta ) - f(\xi )| \leqslant C|\eta |\mathop {\sup
}\limits_{{B_{b|\eta |}}(\xi )} |Xf|,
\end{equation}
for all $f \in {C^1}(G)$ and every  $\xi ,\eta  \in G$.
\end{lemma}
\begin{remark}
The constant $b$ in Lemma 2.8 can be taken 1 when the homogenous
norm $| \cdot |$ is changed by the Carnot-Carath\'{e}odory distance,
see [1] for detail. In the sequel we always suppose $b\geq1$  without loss
of generality.
\end{remark}

\section{Proof of main result}
\textbf{Proof of Theorem 1.1.}  We divide the proof into three
steps.

\textbf{Step 1.} Denote  ${B_k} = {B_{{\rho ^k}}}(0),\rho  =
\frac{1} {2}$. By Lemma 2.1, there exists a solution ${u_k} \in
{C^\infty }({B_k}) \cap C({\bar B_k})$ to the boundary value problem
\[\left\{ {\begin{array}{*{20}{c}}
   {L{u_k} = {f_0} = f(0)}  \\
   {{u_k} = u}  \\

 \end{array} } \right.\begin{array}{*{20}{c}}
   {in\;{B_k},}  \\
   {on\;\partial {B_k}.}  \\

 \end{array} \]
Then  ${v_k} = u - {u_k}$ satisfies the Dirichlet boundary value
problem
\[\left\{ {\begin{array}{*{20}{c}}
   {L{v_k} = f - {f_0}}  \\
   {{v_k} = 0}  \\

 \end{array} } \right.\begin{array}{*{20}{c}}
   {in\;{B_k},}  \\
   {on\;\partial {B_k}.}  \\

 \end{array} \]
By Lemma 2.5, we have
\begin{equation}
\mathop {\sup }\limits_{{B_k}} |{v_k}| \leqslant C{\rho
^{2k}}{\omega _f}({\rho ^k}).
\end{equation}

Since ${w_k} = {u_k} - {u_{k + 1}}$  is  $L$-harmonic in ${B_{k +
2}}$, we have by Lemma 2.2 and (3.1) that
\begin{equation}
\mathop {\sup }\limits_{{B_{k + 2}}} |{X_i}{w_k}| \leqslant C{\rho
^{ - k - 2}}\mathop {\sup }\limits_{{B_{k + 1}}} |{w_k}| \leqslant
C{\rho ^{ - k}}\left(\mathop {\sup }\limits_{{B_{k + 1}}} |{v_k}| +
\mathop {\sup }\limits_{{B_{k + 1}}} |{v_{k + 1}}|\right) \leqslant
C{\rho ^k}{\omega _f}({\rho ^k}).
\end{equation}
and
\begin{equation}
\mathop {\sup }\limits_{{B_{k + 2}}} |{X_i}{X_j}{w_k}| \leqslant
C{\rho ^{ - 2k - 4}}\mathop {\sup }\limits_{{B_{k + 1}}} |{w_k}|
\leqslant C{\rho ^{ - 2k}}\left(\mathop {\sup }\limits_{{B_{k + 1}}}
|{v_k}| + \mathop {\sup }\limits_{{B_{k + 1}}} |{v_{k + 1}}|\right)
\leqslant C{\omega _f}({\rho ^k}).
\end{equation}
Applying Lemma 2.6 to  $u \in {C^2}({B_1}(0))$, it gets a homogenous
polynomial ${P_2}(u,0)$  of degree 2 such that for  $1 \leqslant i,j
\leqslant m$,
$${X_i}{P_2}(u,0)(0) = {X_i}u(0)
$$
and
$${X_i}{X_j}{P_2}(u,0)(0) = {X_i}{X_j}u(0).
$$
By (3.1) and Lemma 2.7, we have
\begin{eqnarray}
\mathop {\sup }\limits_{{B_k}} |{u_k} - {P_2}(u,0)| &\leqslant&
\mathop {\sup }\limits_{{B_k}} |u - {u_k}| + \mathop {\sup
}\limits_{{B_k}} |u - {P_2}(u,0)|\nonumber\\
&\leqslant&C{\omega _f}({\rho ^k}){\rho ^{2k}} + o({\rho ^{2k}})
\leqslant o({\rho ^{2k}}).
\end{eqnarray}

Noting  $L{P_2}(u,0) = Lu(0) = f(0) = L{u_k}$, it sees that  ${u_k} -
{P_2}(u,0)$ is $L$-harmonic, and follows by Lemma 2.2 and (3.4) that
\[\mathop {\sup }\limits_{{B_k}} |{X_i}{u_k} - {X_i}{P_2}(u,0)| \leqslant C{\rho ^{ - k}}o({\rho ^{2k}}) = o({\rho ^k})\]
and
\[\mathop {\sup }\limits_{{B_k}} |{X_i}{X_j}{u_k} - {X_i}{X_j}{P_2}(u,0)| \leqslant C{\rho ^{ - 2k}}o({\rho ^{2k}}) = o(1),\]
hence
\begin{equation}
\mathop {\lim }\limits_{k \to \infty } {X_i}{u_k}(0) =
{X_i}{P_2}(u,0)(0) = {X_i}u(0),
\end{equation}
\begin{equation}
\mathop {\lim }\limits_{k \to \infty } {X_i}{X_j}{u_k}(0) =
{X_i}{X_j}{P_2}(u,0)(0) = {X_i}{X_j}u(0).
\end{equation}

For any point $\xi_{0}$  near the origin satisfying  $|{\xi _0}|
\leqslant 1/4{b^2}$, we have
\begin{eqnarray}
&&|{X_i}{X_j}u({\xi _0}) - {X_i}{X_j}u(0)|
\nonumber\\
&&\leqslant|{X_i}{X_j}u({\xi _0}) - {X_i}{X_j}{u_k}({\xi _0})| +
|{X_i}{X_j}{u_k}({\xi _0}) - {X_i}{X_j}{u_k}(0)| +
|{X_i}{X_j}{u_k}(0) - {X_i}{X_j}u(0)| \nonumber\\
&&: = {I_1} + {I_2} + {I_3}.
\end{eqnarray}

\textbf{Step 2.} We now estimate $I_{1},I_{2}$  and  $I_{3}$,
respectively, to prove (1.2).

To estimate  $I_{3}$, let $k$  satisfy  ${\rho ^{2k + 4}} \leqslant
|{\xi _0}|: = {d_0} \leqslant {\rho ^{2k + 3}}$. It shows by (3.3)
and (3.6) that
\begin{equation}
{I_3} \leqslant \Sigma _{l = k}^\infty |{X_i}{X_j}{u_l}(0) -
{X_i}{X_j}{u_{l + 1}}(0)| \leqslant C\Sigma _{l = k}^\infty
\frac{{{\omega _f}({\rho ^l})}} {{{\rho ^l}}}{\rho ^l} \leqslant
C\int_0^{\sqrt {{d_0}} } {\frac{{\omega (r)}} {r}dr}.
\end{equation}

To estimate  $I_{1}$, we consider the boundary value problem
\[\left\{ {\begin{array}{*{20}{c}}
   {L{u'_k} = {f_{{\xi _0}}} = f({\xi _0})}  \\
   {{u'_k} = u}  \\

 \end{array} } \right.\begin{array}{*{20}{c}}
   {in\;{B_k}({\xi _0}),}  \\
   {on\;\partial {B_k}({\xi _0}).}  \\

 \end{array} \]
Similarly to (3.3) and (3.6), it follows
\begin{equation}
\mathop {\sup }\limits_{{B_{k + 2}}({\xi _0})}
|{X_i}{X_j}{u'_l}({\xi _0}) - {X_i}{X_j}{u'_{l + 1}}({\xi _0})|
\leqslant C{\omega _f}({\rho ^k}),
\end{equation}
\begin{equation}
\mathop {\lim }\limits_{k \to \infty } {X_i}{X_j}{u'_k}({\xi _0}) =
{X_i}{X_j}u({\xi _0}).
\end{equation}
Since   $L({u'_k} - {u_k}) = {f_{{\xi _0}}} - {f_0}$ in  ${B_{k +
2}}({\xi _0})$, it implies
\[L[{u'_k} - {u_k} - \frac{1}
{2}({f_{{\xi _0}}} - {f_0})x_{1,1}^2] = 0,\;in\;{B_{k + 2}}({\xi
_0}).
\]
By Lemma 2.5 and (3.1), we have
\begin{eqnarray}
&&|{X_i}{X_j}{u'_k}({\xi _0}) - {X_i}{X_j}{u_k}({\xi _0})|  \nonumber\\
&&\leqslant |({f_{{\xi _0}}} - {f_0})| +|{X_i}{X_j}{u'_k}({\xi
_0}) - {X_i}{X_j}{u_k}({\xi _0}) - \frac{1}
{2}{X_i}{X_j}({f_{{\xi _0}}} - {f_0})x_{1,1}^2| \nonumber\\
&&\leqslant C{\omega _f}({\rho ^k}) + C{\rho ^{2k}}\mathop {\sup
}\limits_{{B_{k + 2}}({\xi _0})} |{u'_k} - {u_k}| + C\mathop {\sup
}\limits_{{B_{k + 2}}({\xi _0})} |({f_{{\xi _0}}} - {f_0})x_{1,1}^2|
\nonumber\\
&&\leqslant C{\omega _f}({\rho ^k}) + C{\rho ^{2k}}\left(C{\rho
^{2k}}{\omega _f}({\rho ^k}) + \mathop {\sup }\limits_{\partial
{B_{k + 2}}({\xi _0})} |u - {u_k}|\right) + C{\rho ^{2k}}{\omega
_f}({\rho
^k}) \nonumber\\
&&\leqslant C{\omega _f}({\rho ^k}).
\end{eqnarray}
With a similar process to (3.8), one has by (3.9), (3.10) and (3.11)
that
\begin{eqnarray}
{I_1} &\leqslant& |{X_i}{X_j}u({\xi _0}) - {X_i}{X_j}{u'_k}({\xi
_0})| + |{X_i}{X_j}{u'_k}({\xi _0}) - {X_i}{X_j}{u_k}({\xi _0})|
\nonumber\\
&\leqslant& \Sigma _{l = k}^\infty |{X_i}{X_j}{u'_l}({\xi _0}) -
{X_i}{X_j}{u'_{l + 1}}({\xi _0})| + C{\omega _f}({\rho
^k})\nonumber\\
&\leqslant& C\int_0^{\sqrt {{d_0}} } {\frac{{\omega (r)}} {r}dr} .
\end{eqnarray}

Finally, let us estimate  $I_{2}$. Since  ${w_k} \in {C^\infty
}({B_{k + 2}})$, we have by Lemma 2.8 that
\begin{equation}
|{X_i}{X_j}{w_k}({\xi _0}) - {X_i}{X_j}{w_k}(0)| \leqslant
C{d_0}\mathop {\sup }\limits_{\tiny{\begin{array}{*{20}{c}}
   {|\eta | < b|{\xi _0}| < {\rho ^{k + 2}}}  \\
   {l = 1,2, \ldots ,m}  \\

 \end{array} } }|{X_i}{X_j}{X_l}{w_k}(\eta )| \leqslant C{d_0}{\rho ^{ - k}}{\omega _f}({\rho
 ^k}).
\end{equation}
On the other hand, it derives
\begin{eqnarray}
&&|{X_i}{X_j}{u_1}({\xi _0}) - {X_i}{X_j}{u_1}(0)|\nonumber\\
&&\leqslant C{d_0}\mathop {\sup }\limits_{\tiny{\begin{array}{*{20}{c}}
   {|\eta | < b|{\xi _0}| < {\rho ^{k + 2}}}  \\
   {l = 1,2, \ldots ,m}  \\
\end{array} }} |{X_l}{X_i}{X_j}({u_1}(\eta ) - P({u_1},0)(\eta
))|\nonumber\\
&&\leqslant C{d_0}\left(\mathop {\sup }\limits_{{B_1}} |{u_1}| +
\mathop
{\sup }\limits_{{B_1}} |P({u_1},0)|\right)\nonumber\\
&&\leqslant C{d_0}\left(\mathop {\sup }\limits_{{B_1}} |u| +
||f|{|_{{L^\infty }}} + \mathop \Sigma \limits_{0 < |I| \leqslant 2}
\mathop {\sup }\limits_{{B_1}} |{X^I}({u_1} - \frac{1}
{2}{f_0}x_{1,1}^2)| + \mathop \Sigma \limits_{0 < |I| \leqslant 2}
\mathop {\sup }\limits_{{B_1}} |\frac{1}
{2}{X^I}({f_0}x_{1,1}^2)|\right)\nonumber\\
&&\leqslant C{d_0}\left(\mathop {\sup }\limits_{{B_1}} |u| +
||f|{|_{{L^\infty }}}\right).
\end{eqnarray}
Then we get by (3.13) and (3.14) that
\begin{eqnarray}
{I_2} &\leqslant& |{X_i}{X_j}{u_{k - 1}}({\xi _0}) - {X_i}{X_j}{u_{k
- 1}}(0)| + |{X_i}{X_j}{w_{k - 1}}({\xi _0}) - {X_i}{X_j}{w_{k -
1}}(0)|\nonumber\\
&\leqslant& |{X_i}{X_j}{u_1}({\xi _0}) - {X_i}{X_j}{u_1}(0)| +
\Sigma _{l = 1}^{k - 1}|{X_i}{X_j}{w_l}({\xi _0}) -
{X_i}{X_j}{w_l}(0)|\nonumber\\
&\leqslant& C{d_0}\left(\mathop {\sup }\limits_{{B_1}} |u| +
||f|{|_{{L^\infty }}} + \Sigma _{l = 1}^{k - 1}\frac{{{\omega
_f}({\rho ^l})}} {{{\rho ^{2l}}}}{\rho ^l}\right)\nonumber\\
&\leqslant& C{d_0}\left( {\mathop {\sup }\limits_{{B_1}} |u| +
||f|{|_{{L^\infty }}} + \int_{\sqrt {{d_0}} }^1 {\frac{{{\omega
_f}(r)}} {{{r^2}}}dr} } \right).
\end{eqnarray}

Substituting (3.8), (3.12), (3.15) into (3.7), we conclude that for
every $\xi_{0}$  satisfying   ${d_0} = |{\xi _0}| \leqslant
1/4{b^2}$, it holds
\begin{eqnarray}
&&|{X_i}{X_j}u({\xi _0}) - {X_i}{X_j}u(0)|\nonumber\\
&&\leqslant C\left( {{d_0}(\mathop {\sup }\limits_{{B_1}(0)} |u| +
||f|{|_{{L^\infty }}} + \int_{\sqrt {{d_0}} }^1 {\frac{{{\omega
_f}(r)}} {{{r^2}}}dr} ) + \int_0^{\sqrt {{d_0}} } {\frac{{{\omega
_f}(r)}} {r}dr} } \right).
\end{eqnarray}
For any $\xi$  and $\eta$  in  ${B_{1/2}}(0)$, $d = d(\xi ,\eta )$,
let us choose $\xi  = {\xi _1}, \ldots ,{\xi _n} = \eta $ such that
\[d({\xi _i},{\xi _{i + 1}}) = d',d' < {d_0},(n - 1)d' = d,\;for\;1 \leqslant i \leqslant n -
1.
\]
By applying (3.16) to those points, we get (1.2).

\textbf{Step 3.} If $f \in {C^{0,\alpha }}({B_1}(0))$, $\alpha \in
(0,1)$, then
\[|f(\xi ) - f(\eta )| \leqslant {[f]_{{C^{0,\alpha }}}}d{(\xi ,\eta )^\alpha },\]
thus
\[{\omega _f}(r) = \mathop {\sup }\limits_{d(\xi ,\eta ) < r} |f(\xi ) - f(\eta )| \leqslant {[f]_{{C^{0,\alpha }}}}{r^\alpha }.\]
Hence it yields from the right side of (1.2) that
\begin{eqnarray*}
&&d\left( {\mathop {\sup }\limits_{{B_1}(0)} |u| + ||f|{|_{{L^\infty
}}} + \int_{\sqrt d }^1 {\frac{{{\omega _f}(r)}} {{{r^2}}}dr} }
\right) + \int_0^{\sqrt d } {\frac{{{\omega _f}(r)}} {r}dr}
\nonumber\\
&& \leqslant d\left( {\mathop {\sup }\limits_{{B_1}(0)} |u| +
||f|{|_{{L^\infty }}} + {{[f]}_{{C^{0,\alpha }}}}\int_{\sqrt d }^1
{\frac{1} {{{r^{2 - \alpha }}}}dr} } \right) + {[f]_{{C^{0,\alpha
}}}}\int_0^{\sqrt d } {\frac{1} {{{r^{1 - \alpha }}}}dr}
 \nonumber\\
&&\leqslant d\mathop {\sup }\limits_{{B_1}(0)} |u| + \frac{d} {{2 -
\alpha }}{[f]_{{C^{0,\alpha }}}}\left( {\frac{1} {{{{(\sqrt d )}^{1
- \alpha }}}} - 1} \right) + \frac{1} {\alpha }{[f]_{{C^{0,\alpha
}}}}{\left( {\sqrt d } \right)^\alpha }
\nonumber\\
&&\leqslant C{d^{\alpha /2}}\left( {\mathop {\sup
}\limits_{{B_1}(0)} |u| + ||f|{|_{{C^{0,\alpha }}}}} \right).
\end{eqnarray*}
and proves (1.3).

If  $f \in {C^{0,1}}({B_1}(0))$, then
\[{\omega _f}(r) = \mathop {\sup }\limits_{d(\xi ,\eta ) < r} |f(\xi ) - f(\eta )| \leqslant {[f]_{{C^{0,1}}}}r\]
and
\begin{eqnarray*}
&&d\left( {\mathop {\sup }\limits_{{B_1}(0)} |u| + ||f|{|_{{L^\infty
}}} + \int_{\sqrt d }^1 {\frac{{{\omega _f}(r)}} {{{r^2}}}dr} }
\right) + \int_0^{\sqrt d } {\frac{{{\omega _f}(r)}} {r}dr}
\nonumber\\
&&\leqslant d\left( {\mathop {\sup }\limits_{{B_1}(0)} |u| +
||f|{|_{{L^\infty }}} + {{[f]}_{{C^{0,1}}}}\int_{\sqrt d }^1
{\frac{1} {r}dr} } \right) + {[f]_{{C^{0,1}}}}\sqrt d \nonumber\\
&&\leqslant d\left( {\mathop {\sup }\limits_{{B_1}(0)} |u| +
||f|{|_{{L^\infty }}}} \right) + {[f]_{{C^{0,1}}}}\sqrt d \left( {1
+ |\sqrt d \log \sqrt d |} \right)\nonumber\\
&&\leqslant {d^{1/2}}\left( {\mathop {\sup }\limits_{{B_1}(0)} |u| +
||f|{|_{{C^{0,1}}}}\left( {1 + |\sqrt d \log \sqrt d |} \right)}
\right),
\end{eqnarray*}
thus (1.4) is obtained.   $\Box$

%%%%%%%% Bibliography %%%%%%%%%%%%%%%%%%%%%%%%%%%%%%%%%%%%%%%%%%%%%%%%%

\end{document}